\documentclass[11pt]{article}
\usepackage{amssymb,amsfonts,amsmath,latexsym,epsf,tikz,url}

\newtheorem{theorem}{Theorem}[section]

\newtheorem{observation}[theorem]{Observation}
\newtheorem{conjecture}[theorem]{Conjecture}
\newtheorem{corollary}[theorem]{Corollary}
\newtheorem{lemma}[theorem]{Lemma}
\newtheorem{remark}[theorem]{Remark}
\newtheorem{definition}[theorem]{Definition}

\newcommand{\proof}{\noindent{\bf Proof.\ }}
\newcommand{\qed}{\hfill $\square$\medskip}

\begin{document}

\title{Distinguishing number and distinguishing index of natural and fractional powers of graphs}

\author{
Saeid Alikhani  $^{}$\footnote{Corresponding author}
\and
Samaneh Soltani
}

\date{\today}

\maketitle

\begin{center}
Department of Mathematics, Yazd University, 89195-741, Yazd, Iran\\
{\tt alikhani@yazd.ac.ir, s.soltani1979@gmail.com}
\end{center}

\begin{abstract}
The distinguishing number (index) $D(G)$ ($D'(G)$) of a graph $G$ is the least integer $d$
such that $G$ has an vertex labeling (edge labeling)  with $d$ labels  that is preserved only by a trivial
automorphism. For any $n \in \mathbb{N}$, the $n$-subdivision of $G$ is a simple graph $G^{\frac{1}{n}}$ which is constructed by replacing each edge of $G$ with a path of length $n$. 
 The $m^{th}$ power of $G$, is a graph with same set of vertices of $G$ and an edge between two vertices if and only if there is a path of length at most $m$ between them.	
	The fractional power of $G$, denoted by $G^{\frac{m}{n}}$ is  $m^{th}$ power of the $n$-subdivision of $G$ or $n$-subdivision of $m$-th power of $G$. In this paper we study the distinguishing number and distinguishing index of natural and fractional powers of $G$. We show that  the natural powers more than two of a graph  distinguished by three edge labels. Also we show that for a connected graph $G$ of order $n \geqslant 3$ with maximum degree $\Delta (G)$,  $D(G^{\frac{1}{k}})\leqslant min\{s: 2^k+\sum^s_{n=3}n^{k-1}\geqslant \Delta (G)\}$ and for $m\geqslant 3$, $D'(G^{\frac{m}{k}})\leqslant 3$. 
	
\end{abstract}


\medskip
\noindent{\bf AMS Subj.\ Class.:} 05C15, 05E18

\section{Introduction}
Let $G = (V ,E)$ be a simple  graph with $n$ vertices.  We use the standard graph notation (\cite{Sandi}). An
automorphism of $G$ is a permutation $\sigma$ of the vertex set of $G$ with the property that, for any vertices $u$ and $v$, we have
$u\sigma\sim v\sigma$ if and only if $u\sim v$ (note that $v\sigma$ denotes the image of the vertex $v$ under the permutation $\sigma$). 
The set of all automorphisms of $G$, with the operation of composition of permutations, is a permutation group
on $V$ and is denoted by $Aut(G)$.  
A labeling of $G$, $\phi : V \rightarrow \{1, 2, \ldots , r\}$, is  $r$-distinguishing, 
if no non-trivial  automorphism of $G$ preserves all of the vertex labels.
In other words,  $\phi$ is $r$-distinguishing if for every non-trivial $\sigma \in Aut(G)$, there
exists $x$ in $V$ such that $\phi(x) \neq \phi(x\sigma)$. 
The distinguishing number of a graph $G$ has defined by Albertson and Collins \cite{Albert} and  is the minimum number $r$ such that $G$ has a labeling that is $r$-distinguishing.  
Similar to this definition, Kalinkowski and Pil\'sniak \cite{Kali1} have defined the distinguishing index $D'(G)$ of $G$ which is  the least integer $d$
such that $G$ has an edge colouring   with $d$ colours that is preserved only by a trivial
automorphism.  These indices  has developed  and number of papers published on this subject (see, for example \cite{soltani,Klavzar,fish}).

If $x$ and $y$ are two vertices of $G$, then the distance $d(x,y)$ between $x$ and $y$,  is defined as the length of a minimum path connecting $x$ and $y$. The eccentricity of a vertex 
$x$ is $ecc(x)=max\{d(x,u):u\in V(G)\}$ and the radius $r$ and the diameter $d$ of $G$ are defined as the minimum and maximum eccentricity among vertices of $G$, respectively.
  A vertex $u$ of $G$ is called the central vertex if $ecc(u)=r$. The set of all central vertices of $G$, denoted by $Z(G)$, is called the center of $G$. 
For $k\in \mathbb{N}$, the $k$-power of $G$, denoted by $G^k$, is defined on the vertex set $V(G)$ by adding edges joining any two distinct vertices $x$ and $y$ with distance at most $k$ \cite{Agnarsson, Kral}. In other words, $E(G^k)=\{xy : 1 \leqslant d_G(x, y) \leqslant k\}$. Also the $k$-subdivision of $G$, denoted by $G^{\frac{1}{k}}$, is constructed by replacing each edge $v_iv_j$ of $G$ with a path of length $k$, say $P_{v_iv_j}$. These $k$-paths are called superedges and any new vertex is an internal vertex, and is denoted by $w^{\{v_i,v_j\}}_l$ if it belongs to the superedge $P_{v_iv_j}$, $i<j$ and has distance $l$ from the vertex $v_i$, where $l \in \{1, 2, \ldots , k-1\}$.  Note that for $k = 1$, we have $G^{1/1}= G^1 = G$, and if the graph $G$ has $v$ vertices and $e$ edges, then the graph $G^{\frac{1}{k}}$ has $v+(k-1)e$ vertices and $ke$ edges.
The fractional power of $G$, denoted by $G^{\frac{m}{n}}$ is  $m^{th}$ power of the $n$-subdivision of $G$ or $n$-subdivision of $m$-th power of $G$ (\cite{Moharram}). 
Note that the graphs $(G^{\frac{1}{n}})^m$ and $(G^m)^{\frac{1}{n}}$ are different graphs. The fractional power of a graph has introduced by Iradmusa in \cite{Moharram}. He has investigated the chromatic number and clique number of fractional power of graphs. Also, he has studied domination number and independent domination number of fractional powers of graphs (\cite{BIMS}).   Effantin and Kheddouci have studied the $b$-chromatic number of natural  power of graphs and  have obtained  the exact
value for the $b$-chromatic number of power graphs of a path and have determined bounds for the $b$-chromatic number of
power graphs of a cycle (\cite {Effantin}). 
In the study of distinguishing number and distinguishing index of graphs,  this naturally raises the question:  What happens to the distinguishing number and the  distinguishing index, when we consider the natural power and the  fractional power of a graph?  
In this paper we would like to answer to this question. As usual we denote the complete graph, path and cycle of order $n$ by $K_n$, $P_n$ and $C_n$, respectively. Also
$K_{1,n}$ is the star graph with $n+1$ vertices.

\medskip
In the next section, we state some results on the automorphism group and the distinguishing number and index  of natural power of a graph and then compute theses two parameters for 
natural powers of paths and cycles.  We show that  the natural powers more than two of a graph  distinguished by three edge labels. 
In Section 3 and 4,  we study the  distinguishing number and
the distinguishing  index of the fractional powers  of graphs, respectively.


 \section{The distinguishing number and the distinguishing index of the natural powers of a graph}
 
 In this section, we consider the natural powers of a graph and study their distinguishing number and  distinguishing index.  
 We begin with the following lemma which follows from  the definition of the power of  graph.
 
 \begin{lemma}\label{lem2}
 	Let $G$ be a connected graph of order $n$ and diameter $d$. Then 
 	\begin{itemize}
 		\item[(i)] For every natural number   $t\geqslant d$, $G^t=K_{n}$. 
 		\item[(ii)]  {\rm (Theorem $1$ in \cite{Hobbs})} Let $k = mn$, where $m$ and $n$ are positive integers. Then
 		$G^k= (G^m)^n$.
 		\item[(iii)] Let $x$ and $y$ be two vertices of $G$ with distance   $d_G(x,y)=kq+r$ where $0\leqslant r<k$. Then $d_{G^k}(x,y)=q +r$. 
 	\end{itemize}
 \end{lemma}

 \begin{theorem}\label{thm6}
 	Let $G$ be a connected graph with radius $r$. Then 
 \begin{itemize}
 \item[(i)] The automorphism group of $G$, $Aut (G)$, is a subgroup of $Aut (G^k)$ for $k\geqslant 2$. 
 \item[(ii)] The automorphism group of $G^{2t-1}$, $Aut (G^{2t-1})$, is a subgroup of $Aut (G^{2t})$ for $1\leqslant t \leqslant r$. 
 \end{itemize}
 \end{theorem}
 \proof
 (i) Since $Aut (G)$ is a group, it is suffices  to show that $Aut (G) \subseteq Aut (G^k)$. Let $f$ be an automorphism of $G$.   It is clear that  $v_{i}$ and $v_{j}$ are adjacent  in $G^k$ if and only if $d_{G}(v_{i},v_j)\leqslant k$ and this is true  if and only if $d_{G}(f(v_{i}),f(v_{j}))\leqslant k$ and so  $ f(v_{i})$ and $f(v_{j})$ are adjacent in $G^k$. So $f\in Aut (G^k)$, and the result follows.
 
 (ii)  If $t=1$, then the result follows from Part (i). So let $t> 1$ and $f$ be an automorphism of $G^{2t-1}$.   It is clear that  $v_{i}$ and $v_{j}$ are adjacent  in $G^{2t}$ if and only if $d_{G}(v_{i},v_j)\leqslant 2t$ and this is true  if and only if $d_{G^{2t-1}}(v_{i},v_j)\leqslant 2$, and again this is true if and only if  $d_{G^{2t-1}}(f(v_{i}),f(v_{j}))\leqslant 2$ and so  $d_{G}(f(v_{i}),f(v_{j}))\leqslant 2t$ and hence  $ f(v_{i})$ and $f(v_{j})$ are adjacent in $G^{2t}$. So $f\in Aut (G^{2t})$, and the result follows.\qed

 By  Parts (i) and (ii) of Theorem \ref{thm6}, we have the following results which are comparison between the distinguishing number of a graph and the distinguishing number of its natural powers:
 
  \begin{corollary}\label{cor6}
 Let $G$ be a connected graph with radius $r$. Then 
 \begin{itemize}
 \item[(i)] For every $k\geqslant 2$,  $D(G)\leqslant D(G^k)$. 
 \item[(ii)] For every $1\leqslant t \leqslant r$, $D(G^{2t-1})\leqslant D(G^{2t})$.
 \end{itemize}
 \end{corollary}

 \begin{theorem}
Let $G$ be a connected graph of order $n$ with  diameter $d$ and radius $r$. If    $Z(G)=\{x_1,\ldots , x_t\}$, $t\geqslant 1$, is  the center  of $G$, then for $0\leqslant i \leqslant d-r$ we have
\begin{equation*}
D(G^{r+i})\geqslant \big\vert  \{x\in V(G) \vert ~ 0\leqslant d_G (x_j, x)\leqslant i ~\textsl{for some $j = 1, \ldots  , t$}\}\big \vert. 
\end{equation*} 
 \end{theorem}
 \proof
 We first prove the case $i=0$. If $i=0$, then we have 
\begin{equation*}
\big\{x\in V(G) \vert ~ d_G (x_j, x)=0 ~\textsl{for some $j = 1, \ldots  , t$}\big\}=\{x_1,\ldots , x_t\}.
\end{equation*}

By definition of central vertex and power graph, the vertices $x_1,\ldots , x_t$ are the only vertices of $G^r$ such that $deg_{G^r}x_j =n-1$ where $1\leqslant j \leqslant t$. So the maps that fix noncentral vertices and act on central vertices as an permutation of $\mathbb{S}_t$, are automorphisms of $G^r$. Because they preserve the adjacency relation in $G^r$. Thus we should have at least $t$ labels  so that we have  a vertex distinguishing labeling that is not preserved by last automorphisms. Therefore $D(G^r)\geqslant t$.

For $i> 0$, the proof is similar. Indeed, the elements of the set $\{x\in V(G) \vert ~ 0\leqslant d_G (x_j, x)\leqslant i ~\textsl{for some $j = 1, \ldots  , t$}\}$ are the only vertices of $V(G)$ such that their degree  is $n-1$ in $G^r$. \qed

A graph $G$ is Hamiltonian connected, if and only if every two distinct vertices of $G$ are joined by a Hamiltonian path in $G$ (\cite{Chartrand},\cite{ore}). 
The following theorem implies that the cube of every connected graph is Hamiltonian connected  (see also \cite{Karaganis}).

\begin{theorem} {\rm \cite{Sekanina}}\label{thm9} If $G$ is a connected finite graph, then $G^3$ is Hamiltonian connected.
\end{theorem}

We recall that a traceable graph is a graph that possesses a Hamiltonian path.

\begin{theorem}{\rm \cite{nord}}\label{thm10}
	If $G$ is a traceable graph of order $n \geqslant 7$, then $D'(G) \leqslant 2$.
\end{theorem}
The assumption $n \geqslant 7$ is substantial in this  theorem, because  $D'(K_{3,3}) = 3$.
The following corollary shows that the natural powers more than two of a graph of order 
at least seven can be distinguished by  two labels. 

\begin{corollary}\label{cor2}
	If $G$ is a connected finite graph of order $n \geqslant 7$, then for any  $i\geqslant 3$, $D'(G^i) \leqslant 2$.
\end{corollary}
\proof
It can follows from Theorem \ref{thm9} that $G^i$, $i\geqslant 3$, is Hamiltonian connected, and so it is a traceable graph. Since the order of graph is  $n \geqslant 7$, so by Theorem \ref{thm10} we have $D'(G^i) \leqslant 2$.\qed

\begin{remark}\label{rem1} 
	If $G$ is a connected finite graph of order $1 \leqslant n \leqslant 5$ and $G$ is not a path graph, then the diameter of $G$ is less than or equal to three. So we have $D'(G^i)=D'(K_{n}) =3$ for $i\geqslant 3$. For $n=6$, using table of graphs, observe that the diameter of $G$ is less than or equal to $3$ except eight cases. In that eight cases the diameter of $G$ is $4$ ($G$ is not a path graph), and so $D'(G^i)=D'(K_6) =2$ for $i\geqslant 4$. It can be  easily  computed that $D'(G^3)\leqslant 3$ for that eight cases. 
\end{remark}

Motivated by Corollary \ref{cor2} and  Remark \ref{rem1}, we shall prove that the natural powers more than two of any graph  can be distinguished by  three edge labels. To do this, we need to consider the distinguishing index (and distinguishing number) of natural powers of paths.

  \begin{theorem}\label{thm8}
  	Let $n\geqslant 4$ and $k\geqslant 2$ be  integers. The distinguishing number of the path  $P_n^k$ with diameter $d$ and radius $r$ is as follows:
  	
  	\begin{equation*}
  	D(P_n^{k})=\left\{
  	\begin{array}{ll}
  	2& 1\leqslant k \leqslant r,\\
  	2k-n & r+1\leqslant k \leqslant d.
  	\end{array}\right.
  	\end{equation*}
  \end{theorem}
  \proof
  It can be seen that the degree sequence of $P_n^k$ for $1\leqslant k \leqslant r-1$ is as follows (note that the $i$-th term of the degree sequence is degree of the $i$-th vertex of $P_n$ from left side and the vertices of $P_n$ is denoted by $x_1,\ldots , x_{n}$):
  
  \begin{equation*}
  \{deg_{P_n^k}x_i\}_{i=1}^{n}=\{k,k+1,\ldots , k+(k-1), \underbrace{2k,\ldots , 2k}_{(n-2k)-times},k+(k-1),\ldots , k+1,k\}.
  \end{equation*}

  Also the degree sequence of $P_n^r$ is as follows 
  \begin{equation*}
  \{deg_{P_n^r}x_i\}_{i=1}^{n}=\left\{
  \begin{array}{ll}
  \{r,r+1,\ldots, n-1, n,n-1,\ldots, r+1,r\} & \textsl{$n$ is even,}\\
  \{r,r+1,\ldots , n-1, n,n,n-1,\ldots , r+1,r\} & \textsl{$n$ is odd.}
  \end{array}\right.
  \end{equation*}

  For $P_n^k$, $1\leqslant k \leqslant r$, by assigning the two vertices of  degree $i$, $k\leqslant i \leqslant 2k$, the labels $1$ and $2$, and assigning the remaining vertices the label $1$  we have a distinguishing labeling. Because if $f$ is an automorphism of $P^k_n$ such that fixes all the two vertices with  degree $i$, $k\leqslant i \leqslant 2k-1$, then $f$ fixes all vertices. Since we used two  labels, $D(P^k_n)=2$ for $1\leqslant k \leqslant r$. Also, by considering the adjacency relation we have $Aut (P^k_n)=Aut (P_n)$ for $1\leqslant k \leqslant r$. 
  
  The degree sequence of $P^{k}_n$ for $r+1\leqslant k\leqslant d$ is as follows:
  
  \begin{equation*}
  \{deg_{P_n^k}x_i\}_{i=1}^{n}=\{k,k+1,\ldots , n-1, \underbrace{n,\ldots , n}_{(2k-n)-times},n-1,\ldots , k+1,k\}.
  \end{equation*}

  By assigning all the two vertices of  degree $i$, $k\leqslant i\leqslant n-1$, the labels $1$ and $2$, and assigning all vertices of degree $n$, the distinguishing labeling of a complete graph with order of the number of this kind of vertices, we have $D(P_n^{k})=2k-n$ where $r+1\leqslant k \leqslant d$. Also, by considering the adjacency relation and the number of the vertices of induced complete subgraph we can obtain that $\vert Aut(P_n^{k})\vert = 2\vert Aut (K_{2k-n})|$ where $r+1\leqslant k \leqslant d$.\qed

  Using the proof of Theorem \ref{thm8}, we have the following result which is the distinguishing index of natural powers of paths. 
  
  \begin{corollary}\label{cor3}
  	Let $n\geqslant 3$ and $k\geqslant 2$ be integers. The distinguishing number of the path  $P_n^k$  with diameter $d$ and radius $r$ is as follows:
  	
  	\begin{equation*}
  	D'(P_n^{k})=\left\{
  	\begin{array}{ll}
  	D'(P_n)& 1\leqslant k \leqslant r,\\
  	D'(K_{2k-n-2}) & r+1\leqslant k \leqslant d.
  	\end{array}\right.
  	\end{equation*}
  \end{corollary}

Now, we are ready to state the following result which obtain from Corollary \ref{cor2},  Remark \ref{rem1} and Corollary \ref{cor3}. This result implies  that all natural powers more than two of a graph $G$ distinguished by three edge labels.

\begin{corollary}\label{cor4}
	If $G$ is a connected finite graph of order $n$,  then $D'(G^{m})\leqslant 3$ for $m\geqslant 3$.  
\end{corollary}

Corollary \ref{cor4} is true for the  powers more than two of an arbitrary graph. To see what happen to $G^2$, we need some results. The following Lemma is an easy exercise in graph theory literature.

\begin{lemma}\label{lem3} 
	Let $G$ be a connected graph of order $n$. If $|E(G)|\geqslant {n-1\choose 2}+2$, then $G$  has a Hamiltonian  cycle.
\end{lemma}

\begin{theorem}{\rm \cite{Aingworth}}\label{thm11}
	 If $G^2$ is not a complete graph, then the number of edges in which were added to $G$ to construct $G^2$ is at least $n-2$.
	   \end{theorem}

\begin{corollary} Let $G$ be a connected graph of order $n\geqslant 7$ such that  $G^2$ is not a complete graph. If $\vert  E(G)\vert \geqslant \dfrac{1}{2}(n^2+10-5n)$, then $D'(G^2)\leqslant 2$.
	\end{corollary} 

\proof   By Theorem \ref{thm11} we have $\vert E(G^2)\vert \geqslant n-2+ \vert E(G) \vert$. If $n-2+ \vert E(G) \vert \geqslant {n-1\choose 2}+2$ then $G^2$ has a Hamilton cycle by Lemma \ref{lem3}. But $n-2+ \vert E(G) \vert \geqslant  {n-1\choose 2}+2$ concludes that  $\vert  E(G)\vert \geqslant \dfrac{1}{2}(n^2+10-5n)$. Now by Corollary \ref{cor2} we have the result. \qed


 Before ending this section, let to consider the natural  power of a cycle graph of order $n$, i.e., $C_n$, and study its distinguishing number and index. The following theorem is about the automorphism group of powers of a cycle of order $n$, i.e., $C_n$:

 \begin{theorem}\label{thm7}
 	Let $n\geqslant 3$ and $k\geqslant 2$ be integers. Then
 	\begin{equation*}
 Aut (C_n^k)=	\left\{
 	\begin{array}{ll}
  Aut (C_n)&  n > 2k,\\
  Aut (K_n)  & n \leqslant 2k.
 	\end{array}\right.
 	\end{equation*}
 \end{theorem}
 \proof
 Let $V(C_n)=\{x_1,\ldots, x_n\}$ and $E(C_n)=\{\{x_1,x_2\},\{x_2,x_3\},...,\{x_{n},x_1\}\}$.  
 Let  $n> 2k$, we shall  prove that $Aut (C_n^k)= Aut (C^{k-1}_n)$. Let $f\in Aut (C_n^k)$ and $f \notin Aut (C^{k-1}_n)$. Suppose that there exist $x_s,x_t \in V(C_n)$ such that $f(x_s) = x_t$. Since $f \notin Aut (C^{k-1}_n)$, so either  there is $i$, $1\leqslant i \leqslant k-1$ such that $f(x_{s+i})=x_{t+k}$,  or  $f(x_{s-i})=x_{t+k}$. Without loss of generality we can assume that $f(x_{s+i})=x_{t+k}$. 
 Observe that  the distinct vertices 
 \begin{equation*}
 x_{s-(k-i)},x_{s-(k-i-1)},\ldots , x_{s-2},x_{s-1},x_{s+1},x_{s+2},\ldots , x_{s+i-1},x_{s+i+1},\ldots , x_{s+k}
 \end{equation*}
 are  adjacent to both of $x_s$ and $x_{s+i}$ in $C_n^k$. Since the automorphisms preserve the adjacency, so the image of these vertices should be adjacent to both of $x_t$ and $x_{t+k}$. On the other hand $x_{t+1}, x_{t+2},\ldots , x_{t+k-1}$ are the only vertices of $C_n^k$ that are adjacent to both of $x_t$ and $x_{t+k}$. Thus we have
 \begin{equation*}
 \{f(x_{s+i-k}), f(x_{s+i-k+1}),\ldots , f(x_{s}), f(x_{s+1}), \ldots , f(x_{s+k})\}=\{x_t,x_{t+1},\ldots , x_{t+k}\}.
 \end{equation*}
 
 But these two sets have the distinct cardinality, which is a contradiction. So we have $Aut (C_n^k)= Aut (C_n^{k-1})$, and therefore $Aut (C_n^k)=Aut (C_n)$.  
 For the case  $n \leqslant 2k$, we have  $C_n^k= K_n$. Therefore $Aut (C_n^k)= Aut (K_n)$.\qed

 The   following corollary can be verified directly from Theorem \ref{thm7}. 
 \begin{corollary} 	Let $n\geqslant 3$ and $k\geqslant 2$ be integers. We have
 	\begin{enumerate}
 		\item[(i)]
   	$D(C_n^k)=\left\{
 	\begin{array}{ll}
 	D(C_n)&  n > 2k,\\
 	D(K_n)  & n \leqslant 2k.
 	\end{array}\right.$
  	\item[(ii)]  	
 $D'(C_n^k)=\left\{
 	\begin{array}{ll}
 	D'(C_n)&  n > 2k,\\
 	D'(K_n)  & n \leqslant 2k.
 	\end{array}\right.$
 	\end{enumerate}  	
 \end{corollary}

 The characterization of graphs $G$ with $Aut(G)=Aut(G^2)$, is an interesting problem, and we think that most of graphs have this property.  However,
 until now all attempts to characterize these graphs failed, and it remains as open problem. Let to  end this section by posing the following conjecture: 

\medskip
\begin{conjecture}
	\begin{enumerate}
\item[(i)]
If $G$ is a connected graph with diameter $d$ and radius $r$ such that $r< d\leqslant 2r-2$, then $Aut (G)= Aut (G^2)$.
\item[(ii)] 
 If $G$ is a connected bipartite graph with radius $r> 2$, then $Aut (G)= Aut (G^2)$.
\end{enumerate}
\end{conjecture}


\section{Distinguishing number of the fractional power of graphs}

In this section, we study the distinguishing number  of the fractional powers of graphs. 
It can easily be verified that for $n\geqslant 2$ and $k\geqslant 2$, 
$D(P^{\frac{1}{k}}_n)=2$. Also for $n\geqslant 3$ and $k\geqslant 2$, $D(C^{\frac{1}{k}}_n)=2$. 
We state and prove the following lemma to study more on the distinguishing number of the fractional power of graphs.
\begin{lemma}\label{lem1}
Let $G$ be a connected graph of order $n\geqslant 3$ which  is not a cycle. If $f\in Aut (G^{\frac{1}{k}})$,  then $f\vert_{V(G)}\in Aut (G)$.
\end{lemma}
\proof
Let $\{v_1,\ldots , v_n\}$ be the vertex set of $G$. Suppose that $w^{\{v_i,v_j\}}_t$, $1\leqslant t \leqslant k-1$ is an internal vertex of $G^{\frac{1}{k}}$ such that $f(w^{\{v_i,v_j\}}_t)=v_s$  for some $s\in\{1\ldots , n\} $. Since $deg_{G^{\frac{1}{k}}}(w^{\{v_i,v_j\}}_t)=2$, so $deg(v_s)=2$ (note that $deg_G (v_i)=deg_{G^{\frac{1}{k}}}(v_i)$ for all  $i\in\{1\ldots , n\} $.). Let $v_{s_1}$ be the adjacent vertex to $v_s$ in $G$ and $w^{\{v_s,v_{s_1}\}}_1$ be the adjacent vertex to $v_s$ in $G^{\frac{1}{k}}$. Without loss of generality we can assume that 
\begin{equation*}
f(w^{\{v_i,v_j\}}_{t-1})= w^{\{v_s,v_{s_1}\}}_1, f(w^{\{v_i,v_j\}}_{t-2})= w^{\{v_s,v_{s_1}\}}_2 , \ldots , f(w^{\{v_i,v_j\}}_1) =w^{\{v_s,v_{s_1}\}}_{t-1} , f(v_i)= w^{\{v_s,v_{s_1}\}}_t.
\end{equation*}
 
 Then $deg_G(v_i)=2$. Continuing this process, we see that any vertex  of $G$ has degree two, and so $G$ is a cycle, which is a contradiction. \qed
 
 \begin{observation}\label{obser1}
 Let $G$ be a connected graph of order $n\geqslant 3$ which is not a cycle. Let  $i<j$ and $v_i$ and $v_j$ be two adjacent vertices of $G$. Suppose that $f$ is an automorphism of $G^{\frac{1}{k}}$ such that $f(v_i)=v_{i'}$ and $f(v_j)=v_{j'}$. We have two following  cases:
 
 \medskip
 \textbf{Case 1)} If $i' < j'$, then $f(w^{\{v_i,v_j\}}_t)=w^{\{v_{i'},v_{j'}\}}_t$, where $1\leqslant t \leqslant k-1$.
 
  \medskip
 \textbf{Case 2)} If $i' > j'$, then $f(w^{\{v_i,v_j\}}_t)=w^{\{v_{i'},v_{j'}\}}_{k-t}$, where $1\leqslant t \leqslant k-1$.
 \end{observation}
 
\begin{corollary}\label{cor5}
 Let $G$ be a connected graph of order $n\geqslant 3$ which  is not a cycle. Then for every natural number $k$, 
\begin{enumerate}
\item[(i)] $\vert Aut (G^{\frac{1}{k}})\vert = \vert Aut (G) \vert$.
\item[(ii)] $D(G^{\frac{1}{k}})\leqslant D(G)$. 
\end{enumerate}
\end{corollary}
\proof
(i) It follows from Observation \ref{obser1}.

(ii) By Observation \ref{obser1}, if we label the vertices of  the graph $G$  with $D(G)$ labels in a distinguishing way and assign the internal vertices the label $1$, then we have a distinguishing labeling. Therefore $D(G^{\frac{1}{k}})\leqslant D(G)$.  \qed

We need the following definition to obtain more results on the distinguishing number of fractional powers of graphs: 

\begin{definition}{\rm\cite{Kalinowski(2016)}}
The total distinguishing number $D''(G)$ of a graph $G$ is the least number
$d$ such that $G$ has a total colouring with $d$ colours that is preserved only by the identity
automorphism of $G$.
\end{definition}

We also need the following theorem which gives an upper bound for the  total distinguishing number of a graph: 

 \begin{theorem}{\rm\cite{Kalinowski(2016)}}\label{thmkalin}
If $G$ is a connected graph of order $n \geqslant 3$ with maximum degree $\Delta (G)$, then $D''(G) \leqslant \lceil \sqrt{\Delta (G)}\rceil$.
\end{theorem}

 Note that this bound is sharp, because  $D''(K_{1,n})= \lceil \sqrt{n}\rceil$.
 
 \begin{theorem}\label{relation}
 If $G$ is a connected graph of order $n \geqslant 3$, then $D(G^{\frac{1}{2k}})=D''(G^{\frac{1}{k}})$.
 \end{theorem}
 \proof
 It is suffices  to note that $G^{\frac{1}{2k}}$ is constructed by replacing each edge $v_iv_j$ of $G^{\frac{1}{k}}$ with a path of length $2$, say $P_{v_iv_j}$. So if we consider the label of the edges in total labeling of $G^{\frac{1}{k}}$ as the label of internal vertices of  $G^{\frac{1}{2k}}$ then the result follows. \qed

 By Theorems \ref{thmkalin} and \ref{relation} we have the following corollary: 
 \begin{corollary}
 If $G$ is a connected graph of order $n \geqslant 3$ with maximum degree $\Delta (G)$, then $D(G^{\frac{1}{2}})\leqslant \lceil \sqrt{\Delta (G)}\rceil$.
 \end{corollary}
 
Now  we want to obtain a better upper bound for $D(G^{\frac{1}{k}})$. For this purpose, let $S^G_k(x)$, $k\geqslant 0$ denote a sphere of radius $k$ with a center $x$, i.e., the set of all vertices of $G$ at distance $k$ from $x$.
 
 \begin{theorem}\label{thm2}
 If $G$ is a connected graph of order $n \geqslant 3$ with maximum degree $\Delta (G)$, then $D(G^{\frac{1}{k}})\leqslant min\{s: 2^k+\sum^s_{n=3}n^{k-1}\geqslant \Delta (G)\}$.
 \end{theorem}
 \proof
 The basic idea  follows from  the proof of Theorem \ref{thmkalin} (see  \cite{Kalinowski(2016)}).
  If $G$ is  a cycle of order $n$, then since for $n\geqslant 3$ and $k\geqslant 2$, $D(C^{\frac{1}{k}}_n)=2$  we have the result.  Let $v_0$ be a vertex of degree $\Delta (G)$ and $G\neq K_{1,n-1}$ (the case $G = K_{1,n-1}$ has been considered in next theorem). So $S^{G}_2(v_0)$ is nonempty. 
 We label $v_0$ with the label $2$ and consider a BFS tree $T$ of $G$ rooted at $v_0$. We will first label the vertices of the tree $T^{\frac{1}{k}}$. For a given vertex $v$ of $G$, denote $M_t(v)=\{(w^{\{v,u\}}_1,\ldots , w^{\{v,u\}}_{k-1},u) : vu\in E(G)\}$ (note that the vertices $w^{\{v,u\}}_1,\ldots , w^{\{v,u\}}_{k-1}$ are internal vertices on $P_{vu}$). Let $S^{G}_1(v_0)=\{v_1,\ldots , v_p\}$.
  Without loss of generality we can assume that $v_1$ has a neighbour in $S^{G}_2(v_0)$. We label both $k$-ary $(w^{\{v_0,v_1\}}_1,\ldots , w^{\{v_0,v_1\}}_{k-1},v_1)$ and $(w^{\{v_0,v_2\}}_1,\ldots , w^{\{v_0,v_2\}}_{k-1},v_2)$ with a $k$-ary $(1,\ldots , 1)$. Then we label each $k$-ary of $M_t(v_0)\setminus \{(w^{\{v_0,v_1\}}_1,\ldots , w^{\{v_0,v_1\}}_{k-1},v_1),(w^{\{v_0,v_2\}}_1,\ldots , w^{\{v_0,v_2\}}_{k-1},v_2) \}$ with a distinct $k$-ary of labels different from $(1,\ldots , 1)$. 
 Thus $(1,\ldots , 1)$ appears twice as a $k$-ary of labels in $M_t(v_0)$. We will then label the vertices of graph $G$ in such a way that $v_0$ will be the only vertex of $G$ labeled with the label $2$ such that $k$-ary $(1,\ldots , 1)$ appears twice in the  $M_t(v_0)$. Hence  $v_0$ will be fixed by every automorphism of $G^{\frac{1}{k}}$ preserving the labeling. 
 Therefore  all vertices in $S^{G}_1(v_0)$ and all internal vertices on $P_{v_0v_1},\ldots , P_{v_0v_p}$ will also be fixed, except, possibly $(w^{\{v_0,v_1\}}_1,\ldots , w^{\{v_0,v_1\}}_{k-1},v_1)$ and $(w^{\{v_0,v_2\}}_1,\ldots , w^{\{v_0,v_2\}}_{k-1},v_2)$. To distinguish them we label the sets $\{(w^{\{v_1,u\}}_1,\ldots , w^{\{v_1,u\}}_{k-1},u)\in M_t(v_1) : u\in S^{G}_2(v_0)\}$ and $\{(w^{\{v_2,u\}}_1,\ldots , w^{\{v_2,u\}}_{k-1},u)\in M_t(v_2) : v_2u\in E(T), u\in S^{G}_2(v_0)\}$ with two distinct sets of $k$-ary of labels (this is possible because each of these sets contains at most $\Delta (G)-1$ elements, and we have $\Delta (G)$ distinct $k$-ary of labels). 
 Therefore  every internal vertex on the superedges $P_{v_0v_1},\ldots , P_{v_0v_p}$ and $P_{v_1u}, P_{v_2u}$ where $u\in S^{G}_2(v_0)$ will be fixed by every automorphism of $G^{\frac{1}{k}}$ preserving the presented  labeling. For $i=3,\ldots , p$, we then label all elements of $\{(w^{\{v_i,u\}}_1,\ldots , w^{\{v_i,u\}}_{k-1},u) : v_iu\in E(T) , u\in S^{G}_2(v_0)\}$ with distinct $k$-ary of labels different from the $(1,\ldots , 1)$. This is possible. Thus  all other vertices in $S^{G}_2(v_0)$ and all internal vertices on $P_{v_0v_1},\ldots , P_{v_0v_p}$ and $P_{v_1u},\ldots , P_{v_pu}$ where $u\in S^{G}_2(v_0)\setminus S^{G}_1(v_0)$ will be also fixed.

 Then we proceed recursively with respect to the radius $k$ of subsequent sphere $S^{G}_k(v_0)$ according to the ordering of the BFS tree $T$. Suppose all vertices of $S^{G}_i(v_0)=\{u_1,\ldots , u_{l_i}\}$, $i=0,1,\ldots , k$ and  all  internal vertices on $P_{v_av_b}$, where $v_a,v_b\in S^{G}_i(v_0)$, are fixed by every automorphism of $G^{\frac{1}{k}}$ preserving labels. For each subsequent vertex $u_j$, $j=1,\ldots , l_k$ we label every $k$-ary  $(w^{\{u_j,u\}}_1,\ldots , w^{\{u_j,u\}}_{k-1},u)$ where $u$ is a descendent of $u_j$ in $T$, with a distinct $k$-ary of labels except for $(1,\ldots , 1)$.
  This is possible because the number of $k$-ary to be labeled is not greater than the number of admissible $k$-ary of labels. Thus  all neighbours of $u_j$ in $S^{G}_{k+1}(v_0)$ and all internal vertices on the superedges $P_{u_ju}$, where $u$ is a descendent of $u_j$ in $T$, will be also fixed.
 
 Finally, we label all remaining vertex in $V(G^{\frac{1}{k}})\setminus V(T^{\frac{1}{k}})$ with the label $2$. It is easy to see that if $v$ is a vertex labeled with the label $2$ such that the $k$-ary $(1,\ldots , 1)$ appears twice in $M_t(v)$, then $v=v_0$. Hence  all vertices of $G^{\frac{1}{k}}$ are fixed by any automorphism of $G^{\frac{1}{k}}$ preserving this labeling. \qed

 Now we shall show that  the inequality of Theorem \ref{thm2} is sharp.
 
 \begin{theorem}
For  $m\geqslant 3$ and $k\geqslant 2$, $D(K^{\frac{1}{k}}_{1,m})=min\{s: 2^k+\sum^s_{n=3}n^{k-1}\geqslant m\}$.
 \end{theorem}
 \proof
Let to denote the one vertex in a part of $K_{1,m}$ by $v_0$ and the remaining vertices by $v_1,\ldots , v_m$. So the label of $v_0$ can be arbitrary (because it is the only vertex of degree $m$). 
 Using $min\{s: 2^k+\sum^s_{n=3}n^{k-1}\geqslant m\}$ labels, we have at least $m$ different $k$-ary $(c_1,\ldots , c_k)$ of labels. If we assign these $m$ different $k$-ary $(c_1,\ldots , c_k)$ to the vertex set of $m$ superedges that replaced with the $m$ edges $v_0v_1,\ldots , v_0v_m$ (except the vertex $v_0$), then we have distinguishing labeling. 
 If we use less than $min\{s: 2^k+\sum^s_{n=3}n^{k-1}\geqslant m\}$ labels then we have less than $m$ different $k$-ary of labels, so there would exist at least two paths (start with $v_0$ and end with $v_i$ and $v_j$) labeled similarly. Thus we can find a non-trivial automorphism preserving such a labeling.\qed
 
 
 \section{Distinguishing index of the fractional powers of graphs}
 
 In this section, we study the distinguishing index   of the fractional powers  of graphs. It can be easily verified that for $n\geqslant 2$ and $k\geqslant 2$, $D'(P^{\frac{1}{k}}_n)=2$ and for $n\geqslant 3$ and $k\geqslant 2$, $D'(C^{\frac{1}{k}}_n)=2$. We begin with the following theorem. 
 
 \begin{theorem}\label{thm3}
 Let $G$ be a connected graph of order $n\geqslant 3$ such that it is not a cycle. 
 For any  $k\geqslant 2$,  $D(G^{\frac{1}{k+1}})\leqslant D'(G^{\frac{1}{k}})$.
\end{theorem}
 \proof
We define a distinguishing vertex labeling for $G^{\frac{1}{k+1}}$ with $D'(G^{\frac{1}{k}})$ labels. Suppose that $P_{v_iv_j}$ is a superedge that has been replaced with the edge $v_iv_j$ in structure of $G^{\frac{1}{k}}$. If we assign the internal vertices have been lied on corresponding superedge in construction of $G^{\frac{1}{k+1}}$, the label of edges $P_{v_iv_j}$ in $G^{\frac{1}{k}}$ and assign the remaining vertices the label $1$, then by Observation \ref{obser1} we can conclude that the labeling is distinguishing. Since we used $D'(G^{\frac{1}{k}})$ labels, the result follows. \qed

 \begin{theorem}\label{thm4}
  Let $G$ be a connected graph of order $n\geqslant 3$ such that it is not a cycle. Then $D'(G^{\frac{1}{2}})\leqslant \lceil \dfrac{-1+\sqrt{1+8D'(G)}}{2} \rceil$.
 \end{theorem}
 \proof
 First using the label of edges, we partition the edge set of $G$. So we have $D'(G)$ classes, each class contains the edges with similar labels. The elements of $[i]$-th class are denoted by $e_{i1},\ldots , e_{is_i}$, where $1\leqslant i \leqslant D'(G)$ and $\sum_{i=1}^{D'(G)}s_i = \vert E(G) \vert$. We know that each edge of $G$ is replaced with a path of length $2$ in $G^{\frac{1}{2}}$. Let the edge $e_{ij}$ in $G$ be replaced with two edges $e^1_{ij}$ and $e^2_{ij}$ in $G^{\frac{1}{2}}$. For $1\leqslant i \leqslant D'(G)$, we assign the distinct pairs $(c_{i1},c_{i2})$ of labels to all new edges $e^1_{ij}$ and $e^2_{ij}$, where $1\leqslant j \leqslant s_i$ such that 
 \begin{itemize}
 \item[(i)] $c_{i1} \neq c_{i2}$ for $1\leqslant i \leqslant D'(G)$.
 \item[(ii)] $\{c_{i1},c_{i2}\}\neq \{c_{i'1},c_{i'2}\}$ for $1\leqslant i, i' \leqslant D'(G)$ and $i \neq i'$.
 \end{itemize}
 
 By Observation \ref{obser1} this labeling is distinguishing. Since the number of labels that have been used, is $min\{s: \sum_{i=1}^s i \geqslant D'(G)\}$, so $D'(G^{\frac{1}{2}})\leqslant min\{s: \sum_{i=1}^s i \geqslant D'(G)\}$. By an easy computation, we see that
 
 \begin{equation*}
 min\{s: \sum_{i=1}^s i \geqslant D'(G)\}=\lceil \dfrac{-1+\sqrt{1+8D'(G)}}{2} \rceil .
 \end{equation*}
 
 Therefore  we have the result.\qed
 
 \medskip
 Let $(c_1, \ldots , c_k)$ be an $k$-ary of labels such that it is not symmetric, i.e., there exists $i$, $1\leqslant i \leqslant k$ such that $c_i \neq c_{k-i}$. Let the minimum number of labels that have been used in construction of $D'(G)$ numbers of such $k$-ary, is $d_k'$. Then we have the following theorem.

 \begin{theorem}\label{thm5}
  Let $G$ be a connected graph of order $n\geqslant 3$ such that it is not a cycle. Then $D'(G^{\frac{1}{k}})\leqslant d_k'$.
 \end{theorem}
 \proof
 First, we partition the edge set of $G$ using the label of edges. So we have $D'(G)$ classes, each class contains the edges with similar labels. The elements of $[i]$-th class are denoted by $e_{i1},\ldots , e_{is_i}$, where $1\leqslant i \leqslant D'(G)$ and $\sum_{i=1}^{D'(G)}s_i = \vert E(G) \vert$. We know that each edge of $G$ is replaced with a path of length $k$ in $G^{\frac{1}{k}}$. Let the edge $e_{ij}$ in $G$ be replaced with two edges $e^1_{ij},\ldots , e^k_{ij}$ in $G^{\frac{1}{k}}$. For $1\leqslant i \leqslant D'(G)$, we assign the above explained distinct $k$-ary $(c_{i1},\ldots , c_{ik})$ of labels to all new edges $e^1_{ij},\ldots , e^k_{ij}$, where $1\leqslant j \leqslant S_i$.  By Observation \ref{obser1} this labeling is distinguishing. Since the number of labels that have been used, is $d_k'$, so we have the result.\qed

 By Theorems \ref{thm2}, \ref{thm3} and \ref{thm5}, we have two following upper bound for $D(G^{\frac{1}{k}})$:
 
 \begin{equation*}
 \left\{
 \begin{array}{l}
 D(G^{\frac{1}{k}})\leqslant \lambda_1 := d'_{k-1},\\
  D(G^{\frac{1}{k}})\leqslant \lambda_2 := min\{s: 2^k+\sum^s_{n=3}n^{k-1}\geqslant \Delta (G)\}.
 \end{array}\right.
 \end{equation*}
 
 This  raises the question ``which upper bound is better"? It depends to graph.  For instance, the upper bound $\lambda_1$  is better than $\lambda_2$ for fractional powers of star  graphs, $K^{\frac{1}{k}}_{1,m}$ and the situation is different for the one half power  of $F_2=K_1+2K_2$, i.e.,  $F_2^{\frac{1}{2}}$.

 Now by Corollary \ref{cor4} and Theorem \ref{thm5} we have  the following result. 

\begin{corollary}
If $G$ is a connected finite graph of order $n\geqslant 3$,  then $D'(G^{\frac{m}{k}})\leqslant 3$ for $m\geqslant 3$ and $k\geqslant 1$. 
\end{corollary}
\proof
 We know  that $G^{\frac{m}{k}}$ means $(G^{\frac{1}{k}})^m$ or $(G^m)^{\frac{1}{k}}$. In the case $(G^{\frac{1}{k}})^m$ the result follows directly from Corollary \ref{cor4}. For the case $(G^m)^{\frac{1}{k}}$ we have $D'((G^m)^{\frac{1}{k}})\leqslant d'_k$ by Theorem \ref{thm5} and $d'_k$ is the minimum number of labels that have been used in construction of $D'(G^m)$ numbers of symmetric $k$-ary. Since $D'(G^m)\leqslant 3$, so $d'_k\leqslant 3$. Therefore we have $D'((G^m)^{\frac{1}{k}})\leqslant 3$. \qed

The following results follows  easily from Corollaries  \ref{cor6} and \ref{cor5}. 
\begin{corollary}
If $G$ is a connected finite graph of order $n\geqslant 3$,  then for $m\geqslant 3$ and $k\geqslant 1$ we have 
\begin{enumerate}
\item[(i)] $D(G^{\frac{1}{k}})\leqslant D((G^{1/k})^m)$.
\item[(ii)] $D((G^m)^{\frac{1}{k}})\leqslant D(G^m)$.
\end{enumerate}
\end{corollary}



\end{document}